\documentclass[12pt]{amsart}
\usepackage{latexsym}		%
\usepackage{amssymb}
\usepackage{epsfig}		%
\usepackage{amsfonts}

\renewcommand{\baselinestretch}{1.1}		%
\normalsize					%

\newcommand{\IP}{\mathbb{P}}                                     
\newcommand{\IR}{\mathbb{R}}                           
\newcommand{\IC}{\mathbb{C}}
\newcommand{\IZ}{\mathbb{Z}}

\newcommand{\M}{\mathcal{M}}
\newcommand{\cP}{\mathcal{P}}

\newcommand{\cC}{\mathcal{C}}

\newcommand{\cO}{\mathcal{O}}

\newcommand{\G}{\mathcal{G}}

\newcommand{\lb}{\mathfrak{b}}
\newcommand{\lk}{\mathfrak{k}}
\newcommand{\g}{\mathfrak{g}}

\newcommand{\h}{\mathfrak{h}}

\newcommand{\lieu}{\mathfrak{u}}

\newcommand{\Leaf}{\mathcal L}		%
\newcommand{\sleaf}{{\Leaf}}		%
\newcommand{\Lia}{{\text{\rm Lie}}}	%

\newcommand{\lt}{\mathfrak{t}}
\newcommand{\lp}{\mathfrak{p}}

\newcommand{\epf}{\hfill$\square$\\}		%

\newcommand{\wt}{\widetilde}
\newcommand{\wh}{\widehat}
\newcommand{\al}{\alpha}
\newcommand{\be}{\beta}

\newcommand{\Sif}{\Sigma_i(\wh F)}

\newcommand{\Syst}{\text{\rm Syst}}
\newcommand{\Sect}{\text{\rm Sect}}
\newcommand{\Ssect}{\wh {\text{\rm Sect}}}

\newcommand{\Ad}{\text{\rm Ad}}

\newcommand{\pr}{\text{\rm pr}}		%

\newcommand{\sym}{\text{\rm Sym}}
\newcommand{\tr}{\text{\rm Tr}}
\newcommand{\Hom}{\text{\rm Hom}}

\newcommand{\diag}{{\text{\rm diag}}}

\newcommand{\Vect}{{\text{\rm Vect}}}
	
\newcommand{\im}{{\text{\rm Im}}}

\newcommand {\flb}{\lbrack\!\lbrack}
\newcommand {\frb}{\rbrack\!\rbrack}

\theoremstyle{plain}
\newtheorem {hypo}{\bf\hspace{-\parindent}Hypothesis}

\newtheorem {thm}{Theorem}
\newtheorem {prop}[hypo]{Proposition}%

\newtheorem {cor}[hypo]{Corollary}%
\newtheorem {lem}[hypo]{Lemma}%

\theoremstyle{definition}
\newtheorem {defn}[hypo]{Definition}%

\theoremstyle{remark}
\newtheorem {rmk}[hypo]{Remark}%

\def\mapright#1{\smash{
	\mathop{\longrightarrow}\limits^{#1}}}
\def\maprightto#1{\smash{
	\mathop{\longmapsto}\limits^{#1}}}

\def\mapdown#1{\Big\downarrow
	\rlap{$\vcenter{\hbox{$\scriptstyle#1$}}$}}

\def\mapup#1{\Big\uparrow
	\rlap{$\vcenter{\hbox{$\scriptstyle#1$}}$}}

\begin{document}

\date{\today}

\title{Stokes Matrices and Poisson Lie Groups}
\author{P. P. Boalch}
\address{S.I.S.S.A.\\
Via Beirut 2-4\\
34014 Trieste\\
Italy}
\email{boalch@sissa.it}

\maketitle

\renewcommand{\baselinestretch}{1.1}		%
\normalsize

%

%

%

%

\section{Introduction}
The purpose of this paper is to 
point out and then draw some consequences of the
fact that 
the Poisson Lie group $G^*$ dual to $G=GL_n(\IC)$ may
be identified with a certain 
moduli space of meromorphic connections over the unit disc
having an irregular singularity at the origin.
($G^*$ will be fully described in Section \ref{sn: PL gps}.)

The key feature of this point of view is that 
there is a holomorphic map
$$\nu : \g^* \longrightarrow G^*$$
from the dual of the Lie algebra to the group $G^*$,
for each choice of
diagonal matrix $A_0$ with distinct eigenvalues---the `irregular
type'. 
This map is essentially the Riemann-Hilbert map or de$\!$ Rham morphism
for such connections (we will call it the `monodromy map');
it is generically a local analytic isomorphism. 
The main result is: 
\begin{thm} \label{thm: poissonness}
The monodromy map $\nu$ is a 
Poisson map for each choice of irregular type,
where $\g^*$ has its standard complex Poisson structure and 
$G^*$ has its canonical complex Poisson Lie group structure, but
scaled by a factor of $2\pi i$.
\end{thm}

This was conjectured, and proved in the simplest case, 
in \cite{thesis}
based on the observation
that the space
of monodromy/Stokes data  
of such irregular singular 
connections `looks like' the group $G^*$, and that the
symplectic leaves match up.

We will give two  applications.
First, although $\nu$ is neither injective or surjective, upon
restricting to the skew-Hermitian matrices $\lk^*\subset \g^*$ 
it becomes injective, at least when $A_0$ is purely imaginary,
i.e. diagonal skew-Hermitian (both $\lk^*$ and $\g^*$ are identified with
their duals using the trace here). 
We also find that the involution $B\mapsto -B^\dagger$ fixing the
skew-Hermitian matrices corresponds under $\nu$ 
to an involution fixing the Poisson Lie
group $K^*$ dual to the unitary group $K=U(n)$.
This leads to:
\begin{thm}	\label{thm: GW ims}
For each purely imaginary irregular type $A_0$  
the monodromy map restricts to a (real) Poisson diffeomorphism
$\lk^*\cong K^*$ from the dual of the Lie algebra of $K$ to the
dual Poisson Lie group 
(with its standard Poisson structure, scaled by a factor of $\pi$).
\end{thm}

Thus we have a new, direct proof of a
theorem of Ginzburg and Weinstein \cite{GinzW}, that $\lk^*$ and $K^*$ are 
(globally) isomorphic as Poisson manifolds.
Such diffeomorphisms enable one to convert  
Kostant's {\em non-linear} convexity theorem (involving the Iwasawa
projection) 
into Kostant's {\em linear} convexity theorem (which is due to
Schur and Horn in the unitary case, and led to the well-known Atiyah,
Guillemin and Sternberg convexity theorem). 
See \cite{LuR} and Section \ref{sn: K+D} below.
Our approach also gives a new proof of a closely related theorem
of Duistermaat \cite{Dui84}.

Secondly (and this was our original motivation) if we restrict to
skew-symmetric (complex) matrices then the corresponding space of
Stokes data 
naturally appears as a moduli
space of $2$-dimensional topological quantum field theories.
This is due to B. Dubrovin: in \cite{Dub95long} the notion of
a {\em Frobenius manifold} is defined as a geometrical/coordinate-free 
manifestation of the WDVV equations of
Witten-Dijkgraaf-Verlinde-Verlinde
governing deformations of $2$D topological field theories.
One of the main results of \cite{Dub95long} is the identification of
the local moduli of semisimple Frobenius manifolds
with the entries of a Stokes matrix:
an upper triangular matrix $S\in U_+$ with ones on the
diagonal.
An intriguing aspect of \cite{Dub95long} was the explicit formula
for a Poisson bracket on this space of matrices in the three
dimensional case:
\begin{equation}\label{dubs pb intro}
S:=
\left(
\begin{matrix}
1 & x & y \\
0 & 1 & z \\
0 & 0 & 1 
\end{matrix}
\right)
\qquad\quad 
\begin{array}{l}
\{ x, y \} = xy - 2z \\
\{ y, z \} = yz - 2x \\   
\{ z, x \} = zx - 2y.
\end{array}
\end{equation}
This Poisson structure is invariant under a natural braid group action
and has two-dimensional symplectic leaves
parameterised by the values of the Markoff polynomial 
$$x^2+y^2+z^2-xyz.$$
For example, the quantum cohomology of the complex projective plane
$\IP^2(\IC)$ is a $3$-dimensional semisimple Frobenius manifold and
corresponds to the point $S=\left(
\begin{smallmatrix}
1 & 3 & 3 \\
0 & 1 & 3 \\
0 & 0 & 1 
\end{smallmatrix}
\right)$.
(The manifold is just the 
complex cohomology $H^*(\IP^2)$ and the Frobenius structure
comes from the `quantum product', deforming the usual cup product.)
This is an integer solution of the Markoff equation
$x^2+y^2+z^2-xyz=0$
and quite surprisingly it follows that the solution of the WDVV equations
corresponding to the quantum cohomology of $\IP^2$ is not an algebraic
function, from Markoff's proof (in the nineteenth century)
that his equation has an infinite number of integer
solutions (\cite{Dub95long} Appendix F).

Recently M.Ugaglia \cite{Ugag} has extended Dubrovin's formula
to the $n\times n$ case (and found that a constant factor of
$-\frac{\pi i}{2}$ is needed in \eqref{dubs pb intro}).
Our aim here is to relate these Poisson structures to the standard 
Poisson structure on $G^*$:

\begin{thm}	\label{thm: fms+plgps}
The involution of $\g^*$ fixing the skew-symmetric matrices
corresponds under the monodromy map to 
an explicit Poisson involution $i_{G^*}:G^*\to G^*$ having 
fixed point set $U_+$.
The standard ($2\pi i$ scaled) Poisson structure on $G^*$ then induces
the Dubrovin-Ugaglia Poisson structure on the fixed point set $U_+$.
\end{thm}

We note that $U_+$ is not embedded in $G^*$ as a subgroup.
The word `induces' here means the following:
If $S\in U_+\subset G^*$ then the tangent space $T_SG^*$ decomposes
into the $\pm 1$ eigenspaces of the derivative of the involution $i_{G^*}$.
The $+1$ eigenspace is $T_SU_+$ and so there is a projection 
$\pr:TG^*\vert_{U_+}\to TU_+$ along the $-1$ eigenspaces.
The `induced' Poisson bivector on $U_+$ is simply the projection of
the Poisson bivector on $G^*$.

In symplectic terms Theorem \ref{thm: fms+plgps} implies that
symplectic leaves of $U_+$ arise as symplectic submanifolds of
symplectic leaves of $G^*$.

There are other ramifications of the identification of the Poisson Lie group 
$G^*$ as a moduli space of connections that we will postpone.
In particular we plan to elucidate in a future publication
the Poisson braid group action on $G^*$, which arises by virtue of it being
identified with a moduli space of meromorphic connections:
the family of moduli spaces
parameterised by the irregular types $A_0$ has a natural flat 
Ehresmann connection on it (the isomonodromy connection---which can usefully be
thought of as a non-Abelian irregular Gauss-Manin connection \cite{smid}).
The holonomy of 
this Ehresmann connection gives a non-linear Poisson braid group
action on $G^*$.
This action is intimately related to the braid group action on $G^*$ described
explicitly 
by De$\!$ Concini-Kac-Procesi \cite{DKP} in their study of representations of
quantum groups at roots of unity.

The organisation of this paper is as follows.
The next two sections give background material.
Section \ref{sn: PL gps} describes the Poisson Lie groups 
$G^*$ and $K^*$, and Section \ref{sn: monod map} describes the monodromy map,
associating Stokes matrices to an irregular singular connection.
At the end of Section \ref{sn: monod map} we make the basic observation
identifying $G^*$ with a space of meromorphic connections.
Sections \ref{sn: mm is p} and  \ref{sn: GW IMs} then give the proofs of 
Theorems \ref{thm: poissonness} and \ref{thm: GW ims} respectively.
Next Section \ref{sn: K+D} gives some more background material on the
convexity theorems and explains how Duistermaat's theorem arises naturally.
Finally Section \ref{sn: FMs+PLgps} proves Theorem \ref{thm: fms+plgps},
relating Frobenius manifolds to Poisson Lie groups.

Although we work throughout with $G=GL_n(\IC)$, the generalisation to
arbitrary complex reductive groups appears to be straightforward. 
We have not made this generalisation here for fear of obscuring the main
features.

\ 

\renewcommand{\baselinestretch}{1}		%
\small
{\em Acknowledgements.}
The proof of Theorem \ref{thm: poissonness} 
is based on the calculation of Poisson
structures on certain spaces of Stokes matrices due to N.Woodhouse 
\cite{W2000}, to whom I am grateful for sending me \cite{W2000} before
publication.
I would also like to thank B.Dubrovin for advice and encouragement.
A. Weinstein's comment on \cite{GinzW} (Lisbon 1999), that they
``didn't know what the map was'', was also encouraging. 

\renewcommand{\baselinestretch}{1.1}		%
\normalsize

\begin{section}{Poisson Lie Groups} \label{sn: PL gps}

A Poisson Lie group is a Lie group $G$ with a Poisson structure on it
such that the multiplication map $G\times G \to G$ is a Poisson map
(where $G\times G$ is given the product Poisson structure).
This notion was introduced by Drinfel'\!d (see \cite{Drin86});
Poisson Lie groups appear as classical limits of quantum groups.
In other words, one quantises a Poisson Lie group to obtain a quantum
group. 
A remarkable feature is that Poisson Lie groups come in dual pairs:
there is another Poisson Lie group $G^*$ `dual' to any given Poisson Lie
group $G$.
In brief this is because the derivative at the identity of the Poisson
bivector on $G$ is a linear map $\g \to \bigwedge^2\g$, and the dual
of this map is a Lie bracket on $\g^*$. 
The Lie group $G^*$ is defined
as a group with this Lie algebra.
In turn, the Poisson bivector on $G^*$ is determined by requiring its
derivative at the identity to be the dual map of the original Lie
bracket on $\g$; the roles of $G$ and $G^*$ are symmetrical, although
the groups $G$ and $G^*$ are often very different. 

A list of examples appears (in infinitesimal form) in \cite{Drin86}.
Here our main interest is the group $G^*$ dual to $GL_n(\IC)$ with
its standard complex Poisson Lie group structure, so we will
proceed immediately to a description of this case, following
\cite{AlekMalk94, DKP, LuR}.
We will see that the Poisson structure on $G^*$ appears as
a non-linear analogue of the standard linear Poisson structure on 
$\g^*$.

\begin{rmk}
It is relevant to recall that Drinfel'\!d was motivated by Sklyanin's
calculation of the Poisson brackets between matrix entries of a
monodromy matrix $M\in G$ and the observation that this Poisson structure
has the Poisson Lie group property (\cite{Drin86} Remark 5). 
The results here are `dual' to this: a space of Stokes
matrices (i.e. the `monodromy data' of an irregular connection)
will be identified, as a Poisson manifold, with $G^*$. 
\end{rmk}

{\bf The Poisson Lie Group $\bold {G^*}$.}\ \  
Let $B_+,B_-$ be the upper and lower triangular Borel subgroups of
$G:= GL_n(\IC)$, let $U_\pm\subset B_\pm$ be the unipotent subgroups
and $T=B_+\cap B_-\subset G$ the subgroup of diagonal matrices. 
The corresponding Lie algebras will be denoted 
$\lb_+, \lb_-, \lieu_+, \lieu_-, \lt$, all subalgebras of
the $n\times n$ complex matrices $\g=\Lia(GL_n(\IC))$.
The Lie algebra of $G^*$ is defined to be the subalgebra
\begin{equation}  \label{liegstar}
\Lia(G^*) := \{ (X_-,X_+)\in \lb_-\times \lb_+ \ \bigl\vert \ 
		\delta(X_-)+\delta(X_+)=0 \} 
\end{equation}
of the product $\lb_-\times \lb_+$,
where $\delta :\g \to \lt$
takes the diagonal part; $(\delta(X))_{ij}=\delta_{ij}X_{ij}$.
This Lie algebra is identified with the (complex) vector space dual of
$\g$ via the pairing:
\begin{equation}	\label{eqn: cx pairing}
\langle (X_-, X_+), Y \rangle := \tr((X_+ - X_-)Y)
\end{equation}
for any $Y\in \g$.
Thus (\ref{liegstar}) specifies a Lie algebra structure on $\g^*$ and
we define $G^*$ to be the corresponding connected and simply connected
complex Lie group.
Concretely:
\begin{equation}	\label{gstar}
 G^* := \{ (b_-,b_+,\Lambda)\in B_-\times B_+\times\lt \ \bigl\vert \  
\delta(b_-)\delta(b_+)=1, \delta(b_+)=\exp(\pi i \Lambda) \}. 
\end{equation}
It is easily seen that this is an $n^2$ dimensional simply connected
(indeed contractible) subgroup of $B_-\times B_+\times\lt$ (where
$\lt$ is a group under $+$) and has the desired Lie algebra.
(Conventionally \cite{AlekMalk94, DKP}
one omits the $\Lambda$ term appearing in \eqref{gstar} and has $G^*$
non-simply connected; the difference---the choice of
$\Lambda$---is quite trivial, but it is the simply connected group that
arises immediately as a moduli space of meromorphic connections.)

The Poisson bivector on $G^*$ may be defined as follows. 
Consider the map 
\begin{equation}	\label{gstartog}
\pi : G^*\to G; \quad (b_-, b_+, \Lambda)\mapsto b_-^{-1}b_+.
\end{equation}
This is a covering of its image, the `big cell' $G^0\subset G$
containing matrices that admit an `LU factorisation'.

\begin{rmk}	\label{rmk: tau}
If we define, for each $k$, a function $\tau_k:G\to \IC$ taking the
determinant
of the top-left $k\times k$ submatrix of $g\in G$
then note that $\tau_k(b_-^{-1}b_+) = \tau_k(e^{2\pi i \Lambda})$ and
one can prove that  
$G_0 = \{g\in G \ \bigl\vert \ \tau_k(g)\ne 0 \ \forall k\}$.
\end{rmk} 

The conjugation action of $G$ on itself restricts to an infinitesimal
action of $\g$ on $G^0$ (since $G^0$ is open in $G$) and
this lifts canonically along $\pi$ to an infinitesimal
action $\sigma$ of $\g$ on $G^*$ (since $\pi$ is a covering map).
By definition $\sigma:\g \to \Vect(G^*)$ is the Lie algebra
homomorphism taking $X\in\g$ to the corresponding fundamental vector field.
This is the (right) infinitesimal {\em dressing action}.
(It is a `left-action'; the adjective `right' distinguishes $\sigma$
from the {\em left dressing action} which is defined by replacing
$b_-^{-1}b_+$ by $b_+b_-^{-1}$ in \eqref{gstartog}.)

Now, to specify the Poisson bivector $\cP\in\Gamma(\bigwedge^2 TG^*)$ it
is sufficient to give the associated bundle map 
$\cP^\sharp: T^*G^*\to TG^*$ such that
$\cP(\al,\be) = \langle \cP^\sharp(\al), \be \rangle$. 
This is defined simply as the
composition of left multiplication and the right dressing action:
$$\cP_p^\sharp=\sigma_p\circ\varphi;\qquad
T_p^*G^*\mapright{\varphi} T^*_eG^* \cong \g 
\mapright{\sigma_p} T_pG^*$$
where $p\in G^*$ and $\varphi:=l_p^*$ is the dual of the derivative of the map
multiplying on the left by $p$ in $G^*$.
That this does indeed define a Poisson Lie group structure on $G^*$ is
proved in \cite{LuR}. (This is really the complexification
of \cite{LuR} and appears in \cite{AlekMalk94, DKP}---also 
our sign conventions for $\cP^\sharp$ and $\sigma$ are opposite
to \cite{LuR}, however 
these differences cancel out in the definition of $\cP$.)
The same bivector is obtained using right multiplications and the left
dressing action.

\begin{rmk}	\label{rmk: lp analogue}
It is worth noting that the standard Poisson structure on $\g^*$ may be defined
analogously in terms of the coadjoint action and the additive
group structure of $\g^*$.
\end{rmk}

Immediately we can deduce the following well-known fact:
\begin{lem}	\label{lem: Gspleaves}
The symplectic leaves of $G^*$ are the connected components of the 
preimages under $\pi$ of conjugacy classes in $G$.
\end{lem}
\begin{proof}
The tangent space to the symplectic leaf through $p\in G^*$ is the
image of $\cP_p^\sharp: T_p^*G^*\to T_pG^*$, which by definition is the
inverse image under $d\pi$ of the tangent space to the conjugacy class
through $\pi(p)$.
\end{proof}

Another fact that was very motivational is as follows.
Although the infinitesimal dressing actions above 
do not integrate to group actions, the restriction to the diagonal subalgebra
of both the left and right dressing actions 
integrates to the following torus action: 
\begin{equation}	\label{eqn: t action}
t\cdot (b_-,b_+,\Lambda) = (tb_-t^{-1},tb_+t^{-1},\Lambda)
\end{equation}
for any $t\in T$ and $(b_-,b_+,\Lambda)\in G^*$.
Moreover this torus action is Hamiltonian:
\begin{lem}[See also \cite{LuR}]
The map 
$$\mu_T : G^*\longrightarrow \lt^*;\quad (b_-,b_+,\Lambda)\longmapsto 
(2\pi i)\Lambda$$ 
is an equivariant moment map for the torus action \eqref{eqn: t action}.
\end{lem}
\begin{proof}
Choose $X\in \lt$ and let $f:G^*\to \IC; (b_-,b_+,\Lambda)\mapsto 
(2\pi i)\tr(X\Lambda)$ be the $X$ component of $\mu_T$.
Observe that the one-form $df$ on $G^*$ is left-invariant and takes the value
$X\in T^*_eG^*\cong \g$ at $e\in G^*$.
Thus by definition $\cP^\sharp(df) = \sigma(X)$.
This says precisely that $f$ is a Hamiltonian for the vector field $\sigma(X)$
generated by $X$.
\end{proof}
  
\begin{rmk}
1) This lemma will also be an immediate consequence of
Theorem \ref{thm: poissonness} (since $\Lambda$ will be
 essentially the diagonal 
part of a matrix $B\in \g\cong\g^*$ 
and this is a moment map for the coadjoint action
of $T$ on $\g^*$).

2) It is intriguing to observe that the sum of the first $k$ entries of
   $\mu_T$ is a logarithm of the map $\tau_k\circ \pi: G^*\to \IC$,
   where $\pi$ is from \eqref{gstartog} and $\tau_k$ from Remark 
   \ref{rmk: tau}. 
\end{rmk}

Having given the intrinsic formulation of the Poisson bivector on $G^*$,
we now derive some useful formulae. 
Fix $p=(b_-,b_+,\Lambda)\in G^*$.
\begin{lem}	\label{lem: dractionformula}
The right infinitesimal dressing action is given, for any $X\in\g$, by
$$X\maprightto{\sigma_p} (b_-Z_-,b_+Z_+, \dot\Lambda)\in T_pG^*$$
where 
$(Z_-,Z_+)\in \Lia(G^*)$ is determined from $X$ by the equation
$$
b_+Z_+b_+^{-1} - b_-Z_-b_-^{-1} = 
b_+Xb_+^{-1} - b_-Xb_-^{-1}
$$
and $\dot\Lambda = \delta(Z_+)/(\pi i)$.
\end{lem}
\begin{proof}
Immediate upon differentiating the map $\pi:G^*\to G$. 
\end{proof} 

\begin{rmk}
Equivalently, one may readily verify that $Z_\pm$ is given by
$$ Z_\pm = X - b_\pm^{-1}\Ad_p^*(X)b_\pm$$
where $\Ad_p^*(X)\in\g$ is the coadjoint action of $G^*$ on the dual
$\g$ of its Lie algebra.
\end{rmk}

\begin{cor}
The Poisson bivector on $G^*$ is given by 
$$\cP_p(\varphi^{-1}(X),\varphi^{-1}(Y)) = \tr((Z_+-Z_-)Y)$$
where $\varphi=l_p^*: T_p^*G^*\mapright{\cong} \g$ is
the isomorphism coming from left multiplication, $X,Y\in \g$ are
arbitrary and $(Z_-,Z_+)\in \Lia(G^*)$ is determined by $X$ as in
Lemma \ref{lem: dractionformula}. In turn, if $\Leaf\subset G^*$ is a
symplectic leaf and $p\in\Leaf$ then the symplectic structure on $\Leaf$ is
given (in the above notation) by
\begin{equation}	\label{mult KKS}
\omega_\Leaf(\sigma_p(X),\sigma_p(Y)) = \tr((Z_+-Z_-)Y).
\end{equation}
\end{cor}
\begin{proof}
Immediate from the definitions. 
\end{proof}

Formula \eqref{mult KKS} is the $G^*$ analogue of the well-known
Kirillov-Kostant formula for the 
symplectic structure on coadjoint orbits in $\g^*$.  

{\bf The unitary case.}\ \ 
Let $K=U(n)\subset G$ be the group of $n\times n$ unitary matrices.
This is the fixed point set of the involution 
$g\mapsto g^{-\dagger} = (g^\dagger)^{-1}$ of $G$.
Let $\lk=\Lia(K)\subset \g$ denote the set of skew-Hermitian
matrices.

On the Poisson Lie group $G^*$ we are led 
(see Lemma \ref{lem: Herm equivariance}) to consider the involution
\begin{equation}	\label{eqn: hermn G^* involn 1}
(b_-,b_+,\Lambda)\quad\mapsto\quad 
(b_+^{-\dagger},b_-^{-\dagger},-\overline\Lambda).
\end{equation}
Let $K^*$ be the fixed point set of this involution.
Clearly $K^*$ is a subgroup of $G^*$. 
Taking the $B_+$ component projects $K^*$ isomorphically onto
the subgroup
\begin{equation} \label{eqn: usual defn}
\{ b=b_+ \in B_+ \ \bigl\vert \ \text{ the diagonal entries of $b$ are real
and positive } \}
\end{equation}
of $B_+$. (This is the usual definition of $K^*$.)
Thus on the level of Lie algebras
$$\Lia(K^*)\cong 
\{ Z=Z_+ \in \lb_+ \ \bigl\vert \ \text{ the diagonal entries of $Z$ are real }
 \}$$
and we identify $\Lia(K^*)$ with the (real) vector space dual of
$\lk$ by the formula
$$\langle Z,X \rangle = \im \tr(ZX)$$
for any $X\in \lk$. 
(This is half the imaginary part of the restriction of the bilinear
form \eqref{eqn: cx pairing}.)
The right infinitesimal dressing action of $\g$ on $G^*$ restricts to
an action of $\lk$ on $K^*$, and moreover this infinitesimal action 
integrates to a group action; the right dressing action of $K$ on
$K^*$.
Two descriptions of this action are as follows.

1) Observe that the map $\pi:G^*\to G$ restricts to a diffeomorphism
$\pi\vert_{K^*}:K^* \to P; b\mapsto b^\dagger b$ onto the set
$P\subset G$ of positive definite Hermitian matrices.
Then the right dressing action is defined as
$$k\cdot b = \pi\vert_{K^*}^{-1}(kb^\dagger bk^{-1})$$
for any $k\in K$ and $b\in K^*$.

2) Alternatively recall the Iwasawa decomposition of $G$.
This says (rephrasing slightly) that the product map
$K\times K^*\to G; (k,b)\mapsto kb$ is a diffeomorphism.
In particular there is a map $\rho:G\to K^*; g=kb\mapsto b$ taking the
$K^*$ component of $g$.
It easy to see then that the right dressing action is also given by
$$ k\cdot b = \rho(bk^{-1}). $$ 

The standard (real) Poisson Lie group structure on $K^*$ can be
defined as for $G^*$ in terms of left multiplication and the right
dressing action (\cite{LuR} Remark 4.12).
In particular the symplectic leaves are the orbits of the dressing
action which, by 1), are isomorphic to spaces of Hermitian matrices
with fixed positive eigenvalues.
One should note that the symplectic structures on the leaves are {\em not} 
$K$ invariant; rather the dressing actions are Poisson---i.e. such that the
action map $K\times K^*\to K^*$ is a Poisson map, where $K$ has its standard
non-trivial Poisson Lie group structure (\cite{LuR} Remark 4.14).
The basic formulae are as follows.

\begin{lem} \label{lem: unitary formulae}
Let $b$ be a point of $K^*$, 
$\Leaf\subset K^*$ the symplectic leaf containing $b$ and
choose $X,Y\in \lk$.
Then, at $b$ the symplectic form on $\Leaf$ and
Poisson bivector on $K^*$ are 
given by 
$$
\omega_\Leaf(\sigma_b(X),\sigma_b(Y)) = \im\tr(ZY) = 
\cP_{K^*}(\varphi^{-1}(X),\varphi^{-1}(Y))$$
where $\sigma_b:\lk\to T_b\Leaf$ is the right dressing action,
$\varphi:T_b^* K^*\mapright{\cong} T_e^*K^*=\lk$ is induced from left
multiplication by $b$ and $Z:= X-b^{-1}\Ad^*_b(X)b\in \Lia(K^*)$.
\end{lem}

\end{section}

\begin{section}{The Monodromy Map} \label{sn: monod map}

Now we will jump and describe some spaces of meromorphic connections.
Choose a diagonal $n\times n$ matrix $A_0$ with
distinct eigenvalues.
Given a matrix $B\in \g$ we will consider the meromorphic connection
\begin{equation}	\label{conn}
\nabla = d - A; \qquad 
A = \left(\frac{A_0}{z^2} + \frac{B}{z}\right)dz
\end{equation}
on the trivial rank $n$ holomorphic vector bundle over the Riemann sphere.
Thus $\nabla$ has an order two pole at $0$ (irregular singularity) 
and (if $B\ne 0$) a first order pole at $\infty$ (logarithmic singularity).
We will call $A_0$ the `irregular type' of $\nabla$ and once fixed,
the only variable is $B$, which we identify with the element
$\tr(B\,\cdot\,)$ of $\g^*$.

In this section we will define a 
moduli space of meromorphic connections $\M(A_0)$
over the closed unit disc $\Delta\subset\IP^1$ 
having principal parts at $0$ of the form
\eqref{conn}. 
Restricting the connections in \eqref{conn} to the unit disc 
will give a map $\g^* \to\M(A_0)$.
Then $\M(A_0)$ will be identified transcendentally, via the irregular
Riemann-Hilbert correspondence, with a space of monodromy data
$M(A_0)$, containing a pair of Stokes matrices and the so-called
`exponents of formal monodromy'.  
As a manifold there will be a simple identification  $M(A_0)\cong
G^*$ between $M(A_0)$ and the Poisson Lie group $G^*$ defined above.
Thus for each $A_0$ (plus a certain discrete choice---of initial sector
and branch of $\log(z)$) the composition 
$$\g^* \to \M(A_0) \mapright{\text{RH}} M(A_0)\cong G^*$$
defines a holomorphic map $\nu:\g^*\to G^*$; the {\em monodromy map}.
Our aim in this section is to fill in the details of this description.

Suppose $\nabla$ is any  meromorphic connection on a rank $n$ vector
bundle $V$ over the unit disc $\Delta$ with an order two pole at $0$ and no
others. Upon choosing a trivialisation of $V$ we find
\begin{equation}	\label{conn2}
\nabla = d - A; \qquad 
A = \left(\frac{A'_0}{z^2} + \frac{B}{z}\right)dz +\Theta
\end{equation}
for some matrices $A'_0, B\in\g$ and a 
matrix $\Theta$ of holomorphic one-forms on $\Delta$.

A framing of $V$ at $0$ is an isomorphism $g_0:V_0\cong \IC^n$
between the fibre of $V$ at $0$ and $\IC^n$.
We will say a connection with framing $(\nabla,V,g_0)$ has
{\em irregular type} 
$A_0$ if we have $A'_0=A_0$ { in any trivialisation $V\cong
\Delta\times\IC^n$ extending the framing $g_0$}.

\begin{defn}
The moduli space $\M(A_0)$ is the set of isomorphism classes of
triples $(\nabla,V,g_0)$ consisting of a meromorphic connection $\nabla$
on a rank $n$ vector bundle $V\to\Delta$ with just one pole,
of second order at $0$,
together with a framing $g_0$ at $0$ in which
$\nabla$ has irregular type $A_0$.
\end{defn}

Concretely, if $\Syst_\Delta(A_0)$ denotes the infinite dimensional space
of connections \eqref{conn2} on the trivial bundle over $\Delta$
having $A'_0=A_0$, then 
(by choosing arbitrary trivialisations of the bundles $V$ extending their
framings $g_0$)
we obtain an isomorphism
$$ \M(A_0)\cong \Syst_\Delta(A_0)/\G_\Delta$$  
where the gauge group $\G_\Delta$ is the group
of holomorphic maps $g: \Delta\to GL_n(\IC)$ with $g(0)=1$.
We will denote the gauge action by square brackets:
$$ g[\nabla] = d-g[A]; \qquad g[A]= gAg^{-1} + (dg)g^{-1}.$$

The remarkable fact is that we can give an explicit description of
$\M(A_0)$ as a complex manifold in terms of the natural monodromy data
for irregular connections: the Stokes matrices and 
exponents of formal monodromy.

\begin{rmk}
{\em Generically}  
a connection \eqref{conn2} (with $A'_0=A_0$) is gauge equivalent to
a  connection of the form \eqref{conn};
indeed \eqref{conn} is often called the `Birkhoff normal form'.
However not every connection can be reduced to this form and
even if possible, 
the form is not unique: the monodromy map is neither injective or
surjective (see \cite{JLP76} for a detailed analysis in the $n=2$ case).
\end{rmk}

{\bf Stokes Matrices.}\ \ 
Here we mainly follow Balser, Jurkat and Lutz \cite{BJL79} and
Martinet and Ramis \cite{MR91}. 
The presentation and notation is close to \cite{smid}.

Let $Q:= -A_0/z$, so that $dQ=A_0dz/z^2$ and write 
$Q(z)=\diag(q_1,\ldots,q_n)$.

\begin{defn}
1) The {\em anti-Stokes directions} at $0$ associated to $A_0$ are the
directions 
along which $e^{q_i-q_j}$ decays most rapidly as $z\to 0$ for some
$i\ne j$. 
(Equivalently they are the directions between pairs of 
eigenvalues of $A_0$, when plotted in the $z$-plane.)
The number of distinct anti-Stokes directions (clearly even) will be
denoted $2l$.

2) The {\em monodromy manifold} $M(A_0)$ is
$$M(A_0) := U_+\times U_-\times \lt,$$
where, for $(S_+,S_-,\Lambda)\in M(A_0)$,
the matrices $(S_+,S_-)$ will be called {\em Stokes matrices} and  
$\Lambda$ is the {\em permuted exponent of formal monodromy}.

\end{defn}

The aim now is to define a surjective map 
$$\wt\nu:\Syst_\Delta(A_0) \mapright{} M(A_0)$$   
having precisely the $\G_\Delta$ orbits as fibres, thereby inducing an
isomorphism $\M(A_0)\cong M(A_0)$; the (irregular) Riemann-Hilbert
isomorphism. The map $\wt\nu$ will be holomorphic
with respect to any finite number of coefficients of $\nabla$ 
(\cite{BJL79} Remark 1.8).
 
The auxiliary choices needed in order to define $\wt \nu$ are:
1) A choice of initial sector $\Sect_0$ at 
$0$ bounded by two adjacent anti-Stokes
directions and 2) A choice of branch of $\log(z)$ on $\Sect_0$.

Given such a choice of initial sector we will label the anti-Stokes
directions $d_1,d_2,\ldots,d_{2l}$ going in a positive sense and 
starting on the positive edge of $\Sect_0$.
We will write 
$\Sect_i=\Sect(d_i,d_{i+1})$ for the open sector swept out by rays
moving from $d_i$ to $d_{i+1}$ in a positive sense. (Indices are taken
modulo $2l$---so $\Sect_0=\Sect(d_{2l},d_1)$.)

Suppose $\nabla=d-A\in\Syst_\Delta(A_0)$.
The first step in defining $\wt\nu(\nabla)$ is to find a formal
transformation simplifying $\nabla$. Some straightforward
algebra yields: 

\begin{lem}[See \cite{BJL79}]	\label{lem: formal series}
There is a unique formal gauge transformation diagonalising $\nabla$
and removing the holomorphic terms. In other words there is a unique
$\wh F\in G\flb z \frb=GL_n(\IC\flb z\frb)$ 
with $\wh F(0)=1$ such that
$$\wh F\left[ \frac{A_0}{z^2}dz + \frac{\delta(B)}{z}dz \right]
=  \frac{A_0}{z^2}dz  +\frac{B}{z}dz + \Theta$$
as formal series, where $\delta(B)$ is the diagonal part of $B$.
\end{lem}

Thus $\nabla$ is {\em formally} isomorphic to the simple diagonal connection
$\nabla^0:= d-\left(\frac{A_0}{z^2} + \frac{\delta(B)}{z}\right)dz$
(the `formal  normal form of $\nabla$').
Clearly the matrix $z^{\delta(B)}e^{Q}$ is a local fundamental
solution for $\nabla^0$ (i.e. its columns are a basis of solutions).
Thus in turn
$\wh F z^{\delta(B)}e^{Q}$ is a {\em formal} fundamental solution for
$\nabla$. 

The radius of convergence of the series $\wh F$ will in
general be zero however so we do not immediately obtain analytic solutions of
$\nabla$.
The way to proceed is via the following result, which is the outcome of
work of many people (see in the references below).

\begin{thm} \label{thm: sums}
1) On each sector $\Sect_i$
there is a canonical way to choose an  invertible 
$n\times n$ matrix of holomorphic functions $\Sif$ such that 
$\Sif [\nabla^0] =  \nabla$.

2) The matrix of functions $\Sif$ can be analytically continued to the
$i$th `supersector'  
$\Ssect_i:= \Sect\left(d_i-\pi/2,d_{i+1}+\pi/2\right)$ 
and then $\Sif$ is asymptotic to
$\wh F$ at $0$ within $\Ssect_i$.

3) If $g\in G\{z\}$ is a germ of a {convergent} gauge
transformation  and $t\in T$ then 
$$\Sigma_i( g \circ \wh F \circ t^{-1}) = g\circ \Sif \circ t^{-1}. $$
\end{thm}
The point is that on $\Sect_i$ there are generally many
holomorphic isomorphisms between $\nabla^0$ and $\nabla$ which are 
asymptotic to $\wh F$ and one is being chosen in a
canonical way; it is in fact characterised by property 2).
The details of the construction of $\Sif$ will not be needed.
There are basically two equivalent ways to define $\Sif$: 
algorithmic (start with some solution and modify it to obtain the
canonical one---see \cite{BJL79,L-R94}), 
or summation-theoretic
(the series $\wh F$ is `$1$-summable' with sum $\Sif$ on the $i$th sector
---see \cite{BBRS91,MalR92,MR91}; the singular directions of the
summation operator are (contained in) the set of anti-Stokes directions).
The directions which bound the supersectors, where the asymptoticity may be
lost, are often referred to as {\em Stokes directions}.
See for example
\cite{Mal95,Was76} regarding asymptotic expansions on sectors. 

Thus we can immediately construct many fundamental solutions of $\nabla$:

\begin{defn} \label{dfn: ffs}
The {\em canonical fundamental solution} of $\nabla$ on $\Sect_i$ is
$$\Phi_i := \Sif z^{\delta(B)}e^Q$$
where (by convention) the 
branch of $\log(z)$ chosen on $\Sect_0$ is extended to the other
sectors in a {\em negative} sense.  
\end{defn}

The Stokes matrices are essentially 
the transition matrices between the canonical
fundamental solution $\Phi_0$ on $\Sect_0$ and $\Phi_l$ 
on the opposite sector $\Sect_l$,
when they are continued along the two possible paths in the punctured
disk joining these sectors. (In fact
these two Stokes matrices encode {\em all} the 
possible transition matrices between any of the
canonical bases of solutions, although this may not be clear from the
definition below---see \cite{BJL79, smid}.)

Some work is required to get these Stokes 
matrices to be in $U_+, U_-$ however, and the standard method is as follows:

\begin{defn} \label{dfn: P&S}
1) The {\em permutation matrix} $P\in G$ associated to the choice of $\Sect_0$ 
is defined by $(P)_{ij}=\delta_{\pi(i)j}$ where $\pi$ is the
permutation of $\{1,\ldots,n\}$ corresponding to the dominance
ordering of $\{e^{q_1},\ldots,e^{q_n}\}$ along the direction $\theta$
bisecting the sector $\Sect(d_1,d_l)$:
$$ \pi(i) < \pi(j) \qquad \Longleftrightarrow  \qquad e^{q_i}/e^{q_j}\to 0 
\text{ as $z\to 0$ along $\theta$}.  $$

2) The {\em Stokes matrices} $(S_+,S_-)$ of $\nabla$ are the unique matrices
such that:

$\bullet$ If $\Phi_l$ is continued in a positive sense to $\Sect_0$
		then $\Phi_l = \Phi_0\cdot P S_- P^{-1}$, and

$\bullet$ If $\Phi_0$ is continued in a positive sense to $\Sect_l$
		then $\Phi_0 = \Phi_l\cdot P S_+ P^{-1} M_0$,
where $M_0:= \exp(2\pi\sqrt{-1}\delta(B))\in T$ is the 
{\em formal monodromy}
of $\nabla$; it is the actual monodromy of the formal normal form $\nabla^0$. 

3) The {\em exponent of formal monodromy} of $\nabla$ is $\delta(B)$
   and the {\em permuted exponent of formal monodromy} is
$\Lambda := P^{-1}\delta(B) P\in\lt$.
\end{defn}

The crucial fact, motivating the definition of $P$, is:  

\begin{lem}	\label{lem: sms are unipotent}
$S_+\in U_+$ and $S_-\in U_-$. 
\end{lem}
\begin{proof}
Observe $\theta$ and $-\theta$ are the bisecting directions of the two
components of the intersection $\Ssect_0\cap\Ssect_l$ of the supersectors.
(Recall from Theorem \ref{thm: sums}, 
for each $i$, $\Sigma_{i}(\wh F)$ is asymptotic to 
$\wh F$ at $0$  when continued within
$\Ssect_i$.)
Thus 
$z^{\delta(B)} e^Q (PS_-P^{-1}) e^{-Q}z^{-\delta(B)} = 
 \Sigma_{0}(\wh F)^{-1}\Sigma_{l}(\wh F)$
is asymptotic to $1$ within the component of
$\Ssect_0\cap\Ssect_{l}$ containing $-\theta$.
The exponentials dominate so we must have
$(PS_-P^{-1})_{ij} = \delta_{ij}$
unless
$e^{q_i-q_j}\to 0$ as $z\to\ 0$ along $-\theta$.
This says, equivalently, that $S_-\in U_-$. 
The argument for $S_+$ is the same once the change in choice of
$\log(z)$ is accounted for. 
\end{proof}

Thus we have now defined the desired map 
$\wt\nu:\Syst_\Delta(A_0) \to M(A_0)$ taking the Stokes matrices and
(permuted) exponents of formal monodromy.
Part 3) of Theorem \ref{thm: sums} implies that the $\G_\Delta$
orbits are contained in the fibres of $\wt \nu$ so that $\wt \nu$
induces a well-defined map $\M(A_0)\to M(A_0)$. 

\begin{thm}[See \cite{BJL79, BV89}]  \label{thm: basic bijn}
The induced map $\M(A_0)\to M(A_0)$ is {bijective}.
\end{thm}
\begin{proof}
For injectivity, suppose two connections $\nabla_1,\nabla_2\in
\Syst_\Delta(A_0)$ have the same Stokes matrices and exponent of
formal monodromy $\delta(B)$.
Let $\wh F_1, \wh F_2$ be the associated formal isomorphisms 
(with the same normal form $\nabla^0$)
and let $\phi_i=\Sigma_i(\wh F_2)\circ\Sigma_i(\wh F_1)^{-1}$
for $i=0$ and $i=l$.
$\phi_i$ is a holomorphic solution of the connection
$\Hom(\nabla_1,\nabla_2)$ asymptotic to $\wh F_2\circ \wh F_1^{-1}$
at $0$ in $\Ssect_i$.
Since the Stokes matrices are equal,  $\phi_0=\phi_l$ on both components
of the intersection $\Ssect_0\cap\Ssect_l$ and so they fit together to
define an isomorphism $\phi$ between $\nabla_1$ and
$\nabla_2$ on the punctured disc. 
By Riemann's removable singularity theorem it follows that
$\phi$ extends across $0$ (and has {\em Taylor} expansion 
$\wh F_2\circ \wh F_1^{-1}$) and so is the desired element of
$\G_\Delta$. 
Surjectivity follows from a result of Sibuya (see \cite{BJL79} Section
6, \cite{BV89} Section 4), together with the (straightforward to
prove) fact that any meromorphic connection germ is gauge equivalent
to the germ of a meromorphic connection on the unit disc.
\end{proof}

Next we observe how the Stokes matrices encode the local monodromy
conjugacy class, and how they behave under the torus
action changing the framing at $0$:
\begin{lem} \label{lem: ccl+t}
If $[(\nabla,V,g_0)]\in\M(A_0)$ has monodromy data $(S_-,S_+,\Lambda)\in
M(A_0)$ then 

1) The monodromy (in the usual sense) of $\nabla$ 
around a simple positive loop in the punctured disc, is conjugate to
$$S_-S_+e^{2\pi i \Lambda}\in G.$$

2) For any $t\in T$, the framed 
connection $(\nabla,V,t\circ g_0)$ has  
monodromy data $(sS_-s^{-1},$
$sS_+s^{-1},\Lambda)$
where $s:=P^{-1}t P\in T$.
\end{lem}
\begin{proof}
1) When continued in a positive sense, $\Phi_0$ becomes 
$\Phi_lP S_+ P^{-1}M_0$ on $\Sect_l$,  which will become  
$\Phi_0P S_-S_+ P^{-1} M_0 = \Phi_0P S_-S_+ e^{2\pi i \Lambda}P^{-1}$
on continuing around, back to $\Sect_0$.

2) Observe that changing $g_0$ to $t \circ g_0$ corresponds to changing
$\wh F$ to $t\wh F t^{-1}$  and so, by 3) of Theorem \ref{thm: sums},
 the canonical solution 
$\Phi_i$ changes to $t\Phi_i t^{-1}$, whence the result is clear.
\end{proof}

{\bf The Monodromy Map.}\ \ 
Combining the maps above we thus obtain a map $\g^*\to M(A)$, taking a
matrix $B$ to the monodromy data at $0$ 
of the connection $d-(A_0/z^2 + B/z)dz$.
The final step is to identify the monodromy manifold $M(A_0)$ with the
Poisson Lie group $G^*$. 
This is motivated by the following simple observation.
Let $\cO\subset \g^*$ be a generic coadjoint orbit and $\cC\subset G$
the conjugacy class $e^{2\pi i \cO}\subset G$ (where $\cO$ is identified with
an adjoint orbit using the trace).

\begin{lem}
If $B \in \cO$ then $S_-S_+e^{2\pi i \Lambda}\in \cC$,  where
$(S_+,S_-,\Lambda)\in M(A_0)$ is the monodromy data at $0$ of 
the connection $\nabla=d-(A_0/z^2 + B/z)dz$.
\end{lem}
\begin{proof}
By Lemma \ref{lem: ccl+t} the local monodromy of $\nabla$ around zero
is conjugate to $S_-S_+e^{2\pi i \Lambda}$. However $\nabla$ has only
one other pole in $\IP^1$: a first order pole at $\infty$
(logarithmic/regular singularity).
The connection $\nabla$ has residue $B$ at infinity and this implies it
has local monodromy conjugate to $e^{-2\pi i B}$ (see e.g. \cite{Was76}).
Clearly a simple positive loop around $\infty$ is a simple {\em negative}
loop around $0$. 
\end{proof}

Now recall from Lemma \ref{lem: Gspleaves} 
that the symplectic leaves of the Poisson Lie group $G^*$
are obtained by fixing the conjugacy class of $b_-^{-1}b_+$.
Thus we are led to the following: 

\begin{defn}
The isomorphism  $M(A_0)\cong G^*$ is defined to be
$(S_+,S_-,\Lambda)\mapsto (b_-,b_+,\Lambda)$ where
$
b_- = e^{-\pi i \Lambda} S_-^{-1},$ and 
$b_+ = e^{-\pi i \Lambda} S_+e^{2\pi i\Lambda}
$, so that $b_-^{-1}b_+ = S_-S_+\exp(2\pi i\Lambda)$.
\end{defn}

Thus we have completed the last step in the definition of the monodromy map
$\nu$ as the composition 
$$\g^* \to \M(A_0) \mapright{\cong} M(A_0)\cong G^*.$$
(This will be streamlined in Section \ref{sn: mm is p}.)
In summary the above considerations lead us to:
\begin{prop}	\label{prop: summary}
For each choice of irregular type $A_0$, 
initial sector and branch of $\log(z)$,
the monodromy map is a holomorphic map $\nu:\g^*\to G^*$ such that:

1) Any generic symplectic leaf of $\g^*$ maps to a symplectic leaf of $G^*$,

2) If $A_0$ and the initial sector are chosen such that the permutation matrix
   $P=1$, then $\nu$ is $T$-equivariant, where $T$ acts on $\g^*$ by the
   coadjoint action and on $G^*$ by the dressing action, as described
   in \eqref{eqn: t action}.
\end{prop}

\begin{rmk}
In \cite{smid} we generalised the well-known Atiyah-Bott construction of
symplectic structures on moduli spaces of holomorphic connections on compact
Riemann surfaces, to the case of meromorphic connections with arbitrary order
poles. 
(Holomorphic connections correspond to complex representations of the
fundamental group of the surface---one needs to incorporate Stokes data in the
general meromorphic case.)
If the surface has boundary a Poisson structure is obtained on the moduli
space, the symplectic leaves of which are specified by fixing monodromy
conjugacy classes on each boundary component (exactly as in the holomorphic
case).
On specialising to the closed unit disc, this gives another a priori
definition of the Poisson structure on $\M(A_0)$.
One can show (as in \cite{smid}) that the monodromy map is Poisson with
respect to this Poisson structure.
Hence (by Theorem \ref{thm: poissonness} and the above
identification $\M(A_0)\cong G^*$) 
we obtain a gauge theoretic construction of the 
standard Poisson Lie group structure on $G^*$.
\end{rmk}

\end{section}

\begin{section}{The Monodromy Map is Poisson} \label{sn: mm is p}

The main step in the proof of Theorem \ref{thm: poissonness} is to see
that the monodromy map restricts to a symplectic map between generic
symplectic leaves, i.e. that it relates the Kirillov-Kostant
symplectic structure and the analogue \eqref{mult KKS} 
on the leaves in $G^*$. 

Thus choose a generic matrix $J\in \g$ (which at this stage only needs to have
the property that no pair of distinct
eigenvalues differ by an integer) and let $\cO$ be the adjoint orbit
of $J$ (which we identify with a coadjoint orbit using the trace and
so $\cO$ inherits a complex symplectic structure).
Let $\cC=\exp(2\pi i \cO)$ be the corresponding conjugacy class and
let $\Leaf=\pi^{-1}(\cC\cap G^0)\subset G^*$ be the symplectic leaf of
$G^*$ over $\cC$ (more precisely each connected component of $\Leaf$
is a symplectic leaf).
Take the symplectic form on $\Leaf$ to be that given by formula 
\eqref{mult KKS} but divided by $2 \pi i$.

Choose an irregular type $A_0$ and initial sector and branch of
$\log(z)$. (For notational simplicity we will assume that these are
chosen such that the corresponding permutation matrix $P$ is the
identity---the extension to the general case is simple.)

We will now associate the following list of data to a matrix $g\in G$:

$\bullet$ Matrices: $B:=gJg^{-1}\in \cO$ and $\Lambda:=\delta(B)\in \lt$,

$\bullet$ Meromorphic connections on the trivial rank $n$ vector
bundle over $\IP^1$:
$$\nabla = d-A,\quad \nabla^0 = d-A^0,\quad \nabla^\infty =
d-A^\infty$$ 
where
$$ A = \left(\frac{A_0}{z^2}+\frac{B}{z}\right)dz,\quad  
A^0 = \left(\frac{A_0}{z^2}+\frac{\Lambda}{z}\right)dz,\quad
A^\infty = \frac{Jdz}{z}.$$

$\bullet$ Formal series:

$\wh F\in G\flb z\frb$ such that $\wh F[A^0]=A$ and $\wh F(0)=1$ (see
Lemma \ref{lem: formal series}), and similarly:

$\wh H\in G\flb z^{-1}\frb$ 
such that $\wh H[A^\infty]=A$ and $\wh H(\infty)=g$.

$\bullet$ Fundamental solutions of $\nabla$:
$$\text{$\Phi:=\Phi_0$ on $\Sect_0$, \quad 
$\Psi:=\Phi_l e^{\pi i \Lambda}$ on $\Sect_l$
(see Definition \ref{dfn: ffs}) and}$$
$$\text{ $\chi:= Hz^J$ on a neighbourhood of $\infty$ slit along $d_1$.}$$

Here the first anti-Stokes ray $d_1$ is extended to $\infty$ and the
chosen branch of $\log(z)$  on $\Sect_0$ is extended to
$\IP^1\setminus d_1$.
The series $\wh H$ is a formal isomorphism at $z=\infty$ between
$\nabla^\infty$ and $\nabla$ and so is a series solution of
$\Hom(\nabla^\infty,\nabla)$; a connection with a simple
pole at $\infty$. 
This implies $\wh H$ is actually 
convergent and defines a holomorphic map $H:\IP^1\setminus\{0\}\to G$.
(See e.g. \cite{Was76} for the existence, uniqueness and convergence of 
$\wh H$.)
Finally we obtain:

$\bullet$ Monodromy data $(b_+,b_-)\in G^*$ and $C\in G$ relating
these fundamental solutions, as indicated schematically in Figure 1.
(For example the arrow $\chi\mapright{C}\Phi$ means that if $\chi$ is
extended along the arrow then $\chi = \Phi\cdot C$ in the domain of
definition of $\Phi$.)

\setcounter{figure}{0}
\begin{figure}[h] 	
\centerline{\input{smapg.pstex_t}}	
\caption{Configuration in $\IP^1$} \label{smapg fig}
\end{figure}

The fact that a simple positive loop around $0$ is also
a simple negative loop around $\infty$ translates into the important
monodromy relation:
\begin{equation} \label{eqn: monod reln}
b_-^{-1}b_+ = Ce^{2\pi i J}C^{-1}.
\end{equation}

Note that $b_\pm$ only depend on $B$ and not on all of $g$ and that by
definition $\nu(B) = (b_-,b_+,\Lambda)\in \sleaf\subset G^*$.

\begin{prop}	\label{prop: mmap is sp}
The restricted monodromy map $\nu:\cO\to\sleaf$ is a symplectic map.
\end{prop}
\begin{pf}
We will now vary the initial matrix $g$ in the procedure above. 
Note that the fundamental solutions (and therefore all the monodromy
data) will vary holomorphically with $g$ (\cite{BJL79} Remark 1.8).
Choose $X_0,Y_0\in \g$ arbitrarily and suppose we have a two parameter
holomorphic family 
$g(s,t)\in G$ with $\dot g g^{-1} = X_0$ and $g'g^{-1} = Y_0$ at $s=t=0$ (for
example $g=e^{X_0t+Y_0s}g_0$).
Generally we will write $\dot M = \frac{\partial M}{\partial t}$ and
$M' = \frac{\partial M}{\partial s}$ and will exclusively be
interested in the point $s=t=0$; this will be tacitly assumed in
all the expressions below. 

By definition $\dot B=[X_0,B]$ and $B'=[Y_0,B]$ and the Kirillov-Kostant 
symplectic structure
on $\cO$ evaluated on these tangents is:
\begin{equation} \label{eqn: KKstr}
\omega_\cO([X_0,B],[Y_0,B]) = \tr([X_0,Y_0]B).
\end{equation}
On the other side, on the leaf $\sleaf\subset G^*$ we have tangents
$$\nu_*(\dot B) = (\dot b_-, \dot b_+, \dot\Lambda) = 
(b_-Z_-, b_+Z_+, \dot\Lambda)$$
where $(Z_-,Z_+):=(b_-^{-1}\dot b_-,b_+^{-1}\dot b_+)\in\Lia(G^*)$ 
(and similarly for $\nu_*(B')$).
The monodromy relation \eqref{eqn: monod reln} implies that if we define 
$X:=-\dot CC^{-1}\in\g$ then
the value of the fundamental vector field of $X$ under the right 
dressing action
is $\nu_*(\dot B)$, i.e. $\sigma_p(X)=\nu_*(\dot B)$ where 
$p:=(b_-, b_+, \Lambda)\in G^*$. Similarly $\sigma_p(Y)=\nu_*(B')$ with
$Y:=-C'C^{-1}$.
Thus formula \eqref{mult KKS} (after rescaling) says 
\begin{equation}	\label{eqn: mult KKstr2}
\omega_\sleaf(\nu_*(\dot B),\nu_*(B'))= 
\frac{1}{2\pi i}
\tr\left( (b_-^{-1}\dot b_- - b_+^{-1}\dot b_+) C' C^{-1} \right).
\end{equation}
Our task is to show that \eqref{eqn: KKstr} and \eqref{eqn: mult
KKstr2} are equal.
This will be accomplished via the following intermediate expression:

\begin{lem}	 \label{lem: KKintegral}
$$\frac{1}{2\pi i}\oint_{\partial\Delta} 
\tr\left(\dot HH^{-1} \nabla (H'H^{-1})\right) = \tr\left([X_0,Y_0]B\right)$$
where $\Delta\subset\IP^1$ is the unit disc $\{z: |z|\le 1\}$ with its
natural orientation and $\nabla$ acts in the adjoint representation:
$\nabla\left(H'H^{-1}\right) = d\left(H'H^{-1}\right) - [A,H'H^{-1}]$.
\end{lem}
\begin{pf}
Recall $H$ is holomorphic on the opposite hemisphere 
$\Delta^+= \IP^1\setminus \{z: |z|< 1\}$ and that 
$H(w) = g + O(w)$ where $w=z^{-1}$.
A direct calculation then gives that, on $\Delta^+$:
\begin{align*}
\tr\left(\dot HH^{-1} \nabla(H'H^{-1})\right) &= 
\tr\left(\dot HH^{-1}[B,H'H^{-1}]\right)\frac{dw}{w} + O(1)dw \\
&= \tr\left(X_0[B,Y_0]\right)\frac{dw}{w} + O(1)dw.
\end{align*}
The lemma now follows immediately from the residue theorem.  \epf
\end{pf} 

\begin{rmk}
In other words this says that the map 
$\cO\to\wh \cO; B\mapsto \nabla\vert_{\partial \Delta}$ is symplectic,
where $\wh \cO$ is the set of connections on the trivial bundle over the
circle $\partial\Delta$ that have monodromy in the
conjugacy class $\cC$. $\wh \cO$ can be naturally identified with a
coadjoint orbit of the central extension of the loop group of $G$ and
so inherits the Kirillov-Kostant symplectic structure, which is (upto scale):
$$\omega_\alpha(\nabla_{\!\alpha}\phi, \nabla_{\!\alpha}\psi)=
\frac{1}{2\pi i}\oint_{\partial\Delta}
\tr(\phi\nabla_{\!\alpha}\psi).$$
In our situation $\alpha  = A\vert_{\partial\Delta} = 
d\chi\chi^{-1}$ and, since $J$ is fixed,
$\dot \chi\chi^{-1}= \dot H H^{-1}$.
Then it follows that
$\dot \alpha = \nabla(\dot H H^{-1})$ and $\alpha' = \nabla(H'H^{-1})$.
\end{rmk}

The strategy now is to re-evaluate the integral in Lemma 
\ref{lem: KKintegral} in terms of the monodromy data.
First note that the integrand is a holomorphic one-form on 
$\IC^* = \IP^1\setminus \{0,\infty\}$, since both $H$ and $\nabla$ are
holomorphic there.
Thus (by Cauchy's theorem) the value of the integral is independent of
the radius of the circle we integrate around: we will calculate the
limit as the radius tends to zero, capitalising on the fact that we
know the asymptotics at $0$ (in appropriate sectors) of $\Phi$ and $\Psi$.

Divide the circle $\partial \Delta_r$ (bounding the disc
of radius $r$ centred at $0$) into two arcs $a^0_r , a^l_r$ by
breaking it at the points $p_r,q_r$ 
of intersection with the directions $\theta$
and $-\theta$ respectively. (Recall $\theta$ was defined to bisect
$\Sect(d_1,d_l)$.)
$a^0_r$ is an arc in a positive sense from $q_r$ to $p_r$ and
is wholly contained in the supersector $\Ssect_0$ (on which we know
$\Phi \sim \wh Fz^\Lambda e^Q$) and 
$a^l_r$ is a positive arc from $p_r $ to $q_r$ 
contained in $\Ssect_l$ (on which we know
$\Psi \sim \wh Fz^\Lambda e^Qe^{\pi i\Lambda}$). 

Define $\varphi$ to be the holomorphic one-form
$$\varphi:= - \tr\left(\nabla(\dot H H^{-1})H'H^{-1}\right)$$
on $\IP^1\setminus\{0,\infty\}$, so (by Stokes theorem and Leibniz) 
$\frac{1}{2 \pi i} \oint_{\partial \Delta}\varphi$ appears in Lemma 
\ref{lem: KKintegral}.

\begin{lem} \label{lem: phiformulae}
On the supersector $\Ssect_0$ we have
$$\varphi = - df_0 + \varepsilon_0,\quad
f_0 = \tr\left( F_0^{-1}\dot F_0 \Lambda' \log_0 z + 
		\frac{1}{2}\dot\Lambda\Lambda'(\log_0 z) ^2 + 
		\Phi^{-1}\dot\Phi C'C^{-1} \right)$$
where we have continued $F_0:=\Sigma_0(\wh F)$ and $\log_0(z)=\log(z)$
from $\Sect_0$ to $\Ssect_0$, and $\varepsilon_0$ is a one-form such
that $\int_{a_r^0}\varepsilon_0\to 0$ as $r\to 0$.
Similarly on $\Ssect_l$ we have
$$\varphi = - df_l + \varepsilon_l,\  
f_l = \tr\left( F_l^{-1}\dot F_l \Lambda' (\log_l z + \pi i)+ 
		\frac{1}{2}\dot\Lambda\Lambda'(\log_l z+ \pi i)^2 + 
		\Psi^{-1}\dot\Psi (b_-C)'(b_-C)^{-1} \right)$$
where we have continued $F_l:=\Sigma_l(\wh F)$ and $\log_l(z)=\log(z)$
from $\Sect_l$ to $\Ssect_l$, and $\int_{a_r^l}\varepsilon_l\to 0$ as $r\to 0$.
\end{lem}
\begin{pf}
First we recall (see e.g. \cite{Was76}) that if $\epsilon$ is a
holomorphic function on $\Ssect_0$ with asymptotic expansion at $0$
a power series 
$\epsilon\sim\sum_0^\infty a_nz^n$, then $\int_{a_r^0}\epsilon dz \to 0$ as
$r\to 0$.

Now $H'H^{-1}=\chi'\chi^{-1}$ and $\chi=\Phi\cdot C$ (when $\chi$ is
extended along $C$'s arrow) and so 
$$-\varphi = \tr\left(\nabla(\dot\Phi\Phi^{-1})\Phi'\Phi^{-1}\right) +
 \tr\left(\nabla(\dot\Phi\Phi^{-1})\Phi C'C^{-1}\Phi^{-1}\right).$$
The second term is $d\tr(\Phi^{-1}\dot \Phi C'C^{-1})$ by Leibniz
and then a direct calculation substituting the definition $\Phi:= F_0z^\Lambda
e^Q$ into the first term, yields
$$-\varphi= df_0+ \tr\left(
d(F^{-1}\dot F)F^{-1}F' + [F^{-1}F',F^{-1}\dot F]A^0 + 
F^{-1}(F'\dot\Lambda-\dot F\Lambda')\frac{dz}{z}\right)$$
where $F=F_0$. Since $F_0\sim \wh F$ in $\Ssect_0$ it follows that
the long expression here is indeed negligible. This proves the first
statement and the second is analogous. \epf
\end{pf}

Thus $\oint_{\partial\Delta_r}\varphi =
(-f_0+f_l)\bigl\vert^{p_r}_{q_r} +  \epsilon_r$
where $\epsilon_r\to 0$ as $r\to 0$.
If we write $v_r=\log_0(p_r)$ then 
$\log_0(q_r)=\log_l(q_r)=v_r-\pi i$ and $\log_l(p_r)=v_r-2\pi i$ and
we find that
\begin{equation} \label{eqn: penult formula}
\oint_{\partial\Delta_r}\varphi 
= \tr\left(\Phi^{-1}\dot \Phi C'C^{-1} - 
 \Psi^{-1}\dot \Psi (b_-C)'(b_-C)^{-1}\right)\Bigl\vert_{p_r}^{q_r} 
- \pi i(2v_r-\pi i)\tr(\dot\Lambda\Lambda') + \epsilon_r.
\end{equation}

\begin{lem} \label{lem: final formulae}
We have 

$\tr\left(\Phi^{-1}\dot \Phi C'C^{-1} - 
 \Psi^{-1}\dot \Psi (b_-C)'(b_-C)^{-1}\right)(q_r) = 
\tr(b_-^{-1}\dot b_- C' C^{-1}) + \pi i v_r\tr\dot\Lambda\Lambda' +
 \epsilon_r$  and 

$\tr\left(\Phi^{-1}\dot \Phi C'C^{-1} - 
 \Psi^{-1}\dot \Psi (b_-C)'(b_-C)^{-1}\right)(p_r) = 
\tr(b_+^{-1}\dot b_+ C' C^{-1}) - \pi i (v_r-\pi i)\tr\dot\Lambda\Lambda' +
 \epsilon_r$ 

\noindent where each $\epsilon_r\to 0$ as $r\to 0$.
\end{lem}
\begin{pf}
Along $-\theta$ we have $\Phi=\Psi\cdot b_-$. Using this to remove
$\Phi$ from the left-hand side of the first formula and expanding 
$(b_-C)'$ yields:
$$\tr\left(b_-^{-1}\dot b_- C' C^{-1}\right) - 
\tr\left(\Psi^{-1}\dot\Psi b_-'b_-^{-1}\right).$$
To deal with the second term here recall that the diagonal part of
$b_-$ is $e^{-\pi i \Lambda}$ so $b_-'b_-^{-1} = -\pi i \Lambda'+n_-$
for some constant strictly lower triangular matrix $n_-$.
Now $\Psi = F_le^{\Lambda(\pi i + \log_lz)}e^Q$ by definition and so
$$\tr\left(\Psi^{-1}\dot\Psi b_-'b_-^{-1}\right) = 
-\pi i\tr(F_l^{-1}\dot F_l\Lambda') + 
\tr\left(F_l^{-1}\dot F_l z^\Lambda e^{\pi i\Lambda}e^Q n_- 
e^{-Q}e^{-\pi i\Lambda}z^{-\Lambda}\right) -
\pi iv_r \tr(\dot\Lambda\Lambda').$$
The first two terms on the right here tend to zero as $z=q_r\to 0$
along $-\theta$ 
(see Lemma \ref{lem: sms are unipotent} for the second term) and so we
have established the first formula. The second formula arises
similarly (using the fact that $\Phi=\Psi\cdot b_+$ along $\theta$) once
we note that
the monodromy relation \eqref{eqn: monod reln} implies 
$(b_-C)'(b_-C)^{-1}=(b_+C)'(b_+C)^{-1}$. \epf
 \end{pf}

Substituting these into \eqref{eqn: penult formula} we happily find
that the $\tr(\dot\Lambda\Lambda')$ terms cancel, so that
$$
\oint_{\partial\Delta}\varphi \quad =\quad 
\lim_{r\to 0} \oint_{\partial\Delta_r}\varphi \quad =\quad 
\tr\left( (b_-^{-1}\dot b_- - b_+^{-1}\dot b_+)C' C^{-1}\right)
$$
thereby completing the proof of Proposition \ref{prop: mmap is sp}. \epf
\end{pf}
\begin{rmk}
The method of Lemma \ref{lem: inj-etc} below can be used to also show 
that the restricted monodromy map $\nu:\cO\to \Leaf$ is {\em injective}.
\end{rmk}

\begin{pfms}{\em(of Theorem \ref{thm: poissonness}).\ \  }
Let $U\subset \g^*$ be the subset of all matrices having distinct
eigenvalues mod $\IZ$. 
Thus $U$ is a regular Poisson submanifold of $\g^*$ and each of its
symplectic leaves is a coadjoint orbit of the type appearing in
Proposition \ref{prop: mmap is sp}.
It follows then (using local Darboux-Weinstein coordinates for
example) that $\nu\vert_U:U\to G^*$ is a Poisson map (where $G^*$ has
its canonical Poisson structure, as defined in Section \ref{sn: PL
gps}, but multiplied by $2\pi i$).
Now choose arbitrary holomorphic functions $f,g$ on $G^*$ and consider
the holomorphic function
$$\nu^*\{f,g\}_{G^*} - \{\nu^*f,\nu^*g\}_{\g^*}$$
on $\g^*$. We have shown this function vanishes on the dense subset 
$U\subset \g^*$ and so it vanishes everywhere. \epf
\end{pfms}
\end{section}

\begin{section}{Ginzburg-Weinstein Isomorphisms} \label{sn: GW IMs}
In this section we will consider the restriction of the monodromy
map to the skew-Hermitian matrices and prove Theorem \ref{thm: GW
ims}.

We will fix the irregular type $A_0$ to be purely imaginary, so that
there are only two anti-Stokes directions; the two halves of the
imaginary axis. We will take $\Sect_0$ to be the sector containing the
positive real axis $\IR_+$ and use the branch of $\log(z)$ 
which is real on $\IR_+$. Thus, by convention, on $\Sect_1$ (the
opposite sector) $\log(z)$ has imaginary part $-\pi$ on the negative
real axis.

In the previous section we explained how to associate monodromy data
$(b_-,b_+,\Lambda,C)\in G^*\times G$ to a matrix $g\in G$, given a
choice of matrix $J$ which has no distinct eigenvalues
differing by an integer.
In other words we have defined a map 
$$\wh \nu: G\times\g''\to
G^*\times G;\quad \wh\nu(g,J) := (b_-,b_+,\Lambda,C).$$
where $\g'' = \{ J\in \g \ \bigl\vert \ 
\text{ if $p\ne q$ are eigenvalues of $J$ then $p-q\notin\IZ$ } \}$.
Note that the set of skew-Hermitian matrices sits inside $\g''$ and
that $\g''$ is {\em open} in $\g$.

\begin{lem}         \label{lem: Herm equivariance}
The extended monodromy map $\wh \nu$ is equivariant as follows:
$$\wh\nu(g^{-\dagger},-J^\dagger) = 
(b_+^{-\dagger},b_-^{-\dagger},-\overline{\Lambda}, C^{-\dagger})$$
where $\wh\nu(g,J) := (b_-,b_+,\Lambda,C)$.
In particular $\nu(-B^\dagger) =
(b_+^{-\dagger},b_-^{-\dagger},-\overline{\Lambda})$ where
$B=gJg^{-1}$, so that if $B$ is skew-Hermitian then $\nu(B)\in
K^*\subset G^*$.
\end{lem} 
\begin{pf}
Let $i:\IP^1\to\IP^1$ denote complex conjugation; $z\mapsto \overline z$.
If $F:U\to G$ is a smooth map, where $U=i(U)\subset \IP^1$ is an open subset
invariant under $i$, then define 
$F^\iota := i^*(F^{-\dagger}) :U\to G$.
Similarly if $\nabla=d-A$ is a connection on the trivial rank $n$
vector bundle over $U$, define $\nabla^\iota = d- A^\iota$ where
$A^\iota := - i^*(A^\dagger)$.
Since both $i$ and $\dagger$ are anti-holomorphic,
$\iota$ takes holomorphic maps/connections to holomorphic
maps/connections.
(Similarly for formal power series, meromorphic connections etc.)
We will repeatedly use the (easily verified) fact that
$ (F[A])^\iota = F^\iota [A^\iota],$
where the square brackets denote the gauge action.
Note that if 
$A = \left(\frac{A_0}{z^2}+\frac{B}{z}\right)dz$ then 
$A^\iota = \left(\frac{A_0}{z^2}-\frac{B^\dagger}{z}\right)dz.$
It follows then that the list of data associated in Section \ref{sn: mm is p}
to $(g^{-\dagger},-J^{\dagger})$ is (in terms of the corresponding 
data associated to $(g,J)$):
$$(-B^\dagger,-\overline\Lambda,\nabla^\iota,(\nabla^0)^\iota,
(\nabla^\infty)^\iota,
\wh F^\iota, 
\wh H^\iota, \Phi^\iota, \Psi^\iota, \chi^\iota).$$
The only subtlety here involves the fundamental solution 
$\Psi:=\Sigma_1(\wh F)z^\Lambda e^Qe^{\pi i \Lambda}$ on $\Sect_1$.
By definition the new $\Psi$ is 
$\Sigma_1(\wh F^\iota)z^{-\overline\Lambda} e^Qe^{-\pi i \overline\Lambda}$.
To see this is $\Psi^\iota$ we just need to observe that on $\Sect_1$
we have 
$(z^\Lambda)^\iota = z^{-\overline \Lambda}e^{-2\pi i \overline
\Lambda}$.
The lemma now follows immediately.
For example, since $\Phi = \Psi b_+$ on the positive imaginary axis,
we have $\Phi^\iota = \Psi^\iota b_+^\iota$ on the negative imaginary
axis, and $b_+$ is constant so $b_+^\iota  = b_+^{-\dagger}$.\epf
\end{pf}

\begin{rmk}
The involution on the monodromy data is much less
attractive when written in terms of the Stokes matrices; 
one is thus led to believe that $G^*$ is a more natural receptacle.
\end{rmk}

Next we examine the injectivity of $\nu\vert_{\lk^*}$. 

\begin{lem} \label{lem: inj-etc}
1) For each $J\in\g''$ the map 
$\wh\nu_J: G \to G^*\times G;\ g\mapsto \wh\nu (g,J),$
is injective.

2) If $h\in G$ then
$\wh\nu (gh^{-1}, hJh^{-1}) = (b_-,b_+,\Lambda,Ch^{-1})$
where $(b_-,b_+,\Lambda,C):= \wh\nu (g,J)$.

3) $\nu\vert_{\lk^*}:\lk^*\to K^*$ is injective and its derivative is
bijective.
\end{lem}
\begin{pf}
Part 1) is similar to Theorem 
\ref{thm: basic bijn}:
Suppose $\wh \nu_J(g_1)=\wh\nu_J(g_2)$.
We will use subscripts $1,2$ to denote the corresponding auxiliary
data.
Thus $\Phi_1, \Phi_2$ denote the corresponding fundamental solutions on
$\Sect_0$. 
Consider the holomorphic matrix 
$X:=\Phi_1\Phi_2^{-1}$. 
$X$ has asymptotic expansion $\wh F_1\wh F_2^{-1}$ on the supersector
$\Ssect_0$. 
On continuation to $\Sect_1$, we find 
$X:=\Psi_1\Psi_2^{-1}$, since $(b_+)_1=(b_+)_2$.
Similarly $X$ is unchanged on return to $\Sect_0$, and on continuation
to $\infty$ it becomes $\chi_1\chi_2^{-1}$.
Thus $X$ has no monodromy around $0$ and has the same asymptotic expansion
$\wh F_1\wh F_2^{-1}$ on $\Ssect_1$.
Riemann's removable singularity theorem then implies $X$ is
holomorphic across $0$ and across $\infty$ (with Taylor
expansions $\wh F_1\wh F_2^{-1}$ and $\wh H_1\wh H_2^{-1}$ respectively).
Thus $X$ is a matrix of holomorphic functions on $\IP^1$ and so
is constant. 
Its value at $0$ is
$\wh F_1(0)\wh F_2^{-1}(0)=1$ and its value at $\infty$
is $(\wh H_1\wh H_2^{-1})\bigl\vert_{z=\infty}=g_1g_2^{-1}$. 

Part 2) is straightforward. For 3) we argue as follows.
We have a commutative diagram:
\begin{equation}	\label{com diagram}
\begin{array}{ccccccc}
  G\times\g'' & \longrightarrow & G^*\times G\times\g''; &\qquad
& (g,J) & \mapsto & (b_-,b_+,\Lambda,C,J) \\
\mapdown{} && \mapdown{}& \qquad
&  \mapdown{} && \mapdown{}\\
\g'' & \longrightarrow & G^*\times \g''; &\qquad 
& gJg^{-1} & \mapsto & (b_-,b_+,\Lambda,CJC^{-1})
\end{array}
\end{equation}
where $(b_-,b_+,\Lambda,C):=\wh\nu(g,J)$.
The top map is injective by 1).
Also 2) implies that the top map takes fibres of the left map into fibres
of the right map (so the bottom map is well-defined) and moreover
distinct fibres go to distinct fibres (so the bottom map is injective).
Now if $B=gJg^{-1}$ is skew-Hermitian then so is $R:=CJC^{-1}$ by 
Lemma \ref{lem: Herm equivariance}.
The monodromy relation \eqref{eqn: monod reln} says
$b_-^{-1}b_+ = e^{2\pi i R}$, so we see $R$ is determined by
$(b_-,b_+)$; the unique Hermitian logarithm of $b_-^{-1}b_+$ is $2\pi iR$.
Thus restricting the bottom map to $\lk^*\subset \g''$ and
`forgetting $R$' on the right-hand side yields an injective map
$\lk^*\to G^*$. This is of course the restriction of the monodromy map
to $\lk^*$.

Finally we must show that $\nu\vert_{\lk^*}$ is bijective on tangent vectors.
First observe the above argument extends to show that there is an open subset
$U\subset \g''$ which contains $\lk^*$ and on which the monodromy map is
injective. (The unique choice of logarithm on the Hermitian matrices extends
uniquely to matrices sufficiently close to being Hermitian.)
Thus $\nu\vert_U:U\to G^*$ is an injective holomorphic map between
equi-dimensional complex manifolds.
This implies it is biholomorphic onto its image (see e.g. \cite{Ran86} Theorem
2.14). It follows immediately that $d\nu\vert_{\lk^*}$ is bijective.
\epf
\end{pf}

All that is left is to consider surjectivity and the Poisson structures:

\begin{pfms}{\em(of Theorem \ref{thm: GW ims}).\ \  }
Choose any point $b\in K^*$ and let $\Leaf\subset K^*$
be its symplectic leaf (dressing orbit of $K$).
The fact that we can diagonalise the Hermitian matrix
$b^\dagger b$ implies we can choose a diagonal element of $\Leaf$,
which we will write as $e^{\pi i J}$ (with $J$ diagonal and purely
imaginary). 
It is straightforward to see that 
$\nu\vert_{\lk^*}(J) = e^{\pi i J}\in K^*$.

Now let $\cO\subset \lk^*$ be the $K$-coadjoint orbit of $J$.
By the monodromy relation \eqref{eqn: monod reln} and Lemma 
\ref{lem: Herm equivariance} we deduce $\nu(\cO)\subset \Leaf$.
Thus $\nu\vert_\cO$ is a smooth map between two equi-dimensional
compact connected manifolds.
As such it has a well-defined degree which, since it is injective, is
$\pm 1$.
This implies $\nu(\cO)= \Leaf$ since non-surjective maps have degree
zero (c.f. \cite{Bott+Tu}).
Thus the monodromy map $\nu$ maps $\lk^*$ {\em onto} $K^*$, as $b$
was arbitrary.

Now we will examine the symplectic structures of $\cO$ and $\Leaf$.
Let $B\in\lk^*$ be the skew-Hermitian matrix with $\nu(B)=b$ and
choose
$g\in K$ such that $B=gJg^{-1}$.
Thus given arbitrary $X_0,Y_0\in \lk$, Proposition \ref{prop: mmap is sp}
says that 
$$\tr([X_0,Y_0]B) = 
\frac{1}{2\pi i} \tr((Z_+-Z_-)Y) = \frac{1}{\pi} \im\tr(ZY)$$
where $Z=Z_+$ and the rest of the notation is as in 
Proposition \ref{prop: mmap is sp}.
In other words (recalling Lemma \ref{lem: unitary formulae}) 
$\nu\vert_\cO : \cO \to \Leaf$ is symplectic (if we divide the
symplectic form on $\Leaf$ by $\pi$).
Arguing as for Theorem \ref{thm: poissonness} it follows that
$\nu\vert_{\lk^*}$ is Poisson (provided we multiply the Poisson
structure on $K^*$ by $\pi$).
Finally since we have also proved the 
derivative of $\nu\vert_{\lk^*}:\lk^*\to K^*$ 
is bijective we deduce this map is indeed a Poisson
diffeomorphism.
\epf
\end{pfms}

\begin{rmk}
The behaviour of the coadjoint action of the maximal diagonal torus $T_K$ of
$K$ under the restricted monodromy map $\nu\vert_{\lk^*}:\lk^*\to K^*$ 
follows from Lemma \ref{lem: ccl+t} and Proposition \ref{prop: summary}:
If the irregular type $A_0$ is such that the permutation matrix $P=1$ then
$T_K$ acts on $K^*$ by the dressing action \eqref{eqn: t action}. In
general one needs to permute this action as indicated in  
Lemma \ref{lem: ccl+t}.
\end{rmk}

\end{section}

\begin{section}{
Kostant's Non-linear Convexity Theorem and the Theorem of Duistermaat}
\label{sn: K+D}

Let $\lp$ be the set of $n\times n$ Hermitian matrices and $P=\exp(\lp)$ the
set of positive definite Hermitian matrices.
Multiplying by $\sqrt{-1}$ identifies $\lp$ with the skew-Hermitian matrices
$\lk$. In turn $\lk\cong\lk^*$ via the trace, so $\lp$ inherits the standard
Poisson structure from $\lk^*$.
The symplectic leaves $\cO\subset \lp$ are the $\Ad(K)$ orbits, consisting of
matrices with the same $n$-tuple of eigenvalues. The map taking the
diagonal part 
$$\delta: \lp\longrightarrow \IR^n$$ 
is a moment map for the adjoint action of the (maximal) diagonal torus $T_K$ 
of $K$.

Schur and Horn proved classically that the set of diagonal entries appearing
in a fixed orbit $\cO$ is a convex polytope; $\delta(\cO)$ is 
the convex hull of the $\sym_n$
orbit of the $n$-tuple of eigenvalues of $\cO$.
Kostant \cite{Kost73} extended this to arbitrary semisimple groups.
Subsequently Atiyah, Guillemin and Sternberg \cite{Ati82, GS82}
put these results into the very
general context of convexity of the images of moment maps for Hamiltonian
torus actions on compact symplectic manifolds.

Now for the non-linear version: Let $\cC=\exp(\cO)\subset P$ be a set of
positive definite Hermitian matrices with fixed eigenvalues.
The {\em Iwasawa projection} $\wh\delta:G\to \IR^n$ is the map 
$$ g=kan\longmapsto \log(a),$$
where $kan$ is the Iwasawa (Gram-Schmidt) decomposition of $g\in G=GL_n(\IC)$
into the product of a unitary matrix $k$ and diagonal positive real matrix $a$
and a unipotent upper-triangular complex matrix $n$.
Clearly $\cC\subset G$.
Kostant's non-linear convexity theorem \cite{Kost73}
says that on $\cC$ the image of the (non-linear)
Iwasawa projection is the same as the convex polytope appearing above: 
$\wh\delta(\cC) = \delta(\cO)$.

What one would like to have is a map $\eta:\lp\to P$ taking each orbit $\cO$
to $\cC=\exp(\cO)$ and converting $\delta$ into $\wh \delta$---i.e. 
such that following diagram commutes:
\begin{equation}	\label{cd: texp}
\begin{array}{ccc}
  \lp & \mapright{\delta} & \IR^n \\
\mapdown{\eta} && \!\bigl\vert\!\bigl\vert \\
  P & \mapright{\wh \delta} & \IR^n.  
\end{array}
\end{equation}
Clearly taking $\eta(X)=e^X$ maps the orbits correctly, but then the diagram
does not commute. 
However one may `twist' the exponential map appropriately:

\begin{thm}[Duistermaat \cite{Dui84}] \label{thm: duistermaat} \ 

\noindent
There is a real analytic map $\psi:\lp\to K$ such that for each $X\in\lp$:

1) $\wh\delta\bigl( \psi(X)^{-1}\cdot\exp(X)\cdot \psi(X) \bigr) 
= \delta(X)$, and

2) The map $\phi_X: k\mapsto k\cdot \psi(k^{-1}Xk)$ is a diffeomorphism from
   $K$ onto $K$.
\end{thm}

Duistermaat's motivation was to reparameterise certain integrals over $K$,
converting terms involving $\wh \delta$ into terms involving 
the linear map $\delta$.
The proof of the existence of such maps $\psi$ in \cite{Dui84} 
is for connected real semisimple groups $G$ (with
finite centre) and involves an indirect homotopy argument.
Our work in the previous sections immediately gives a new proof (in the
case $G=GL_n(\IC)$); one may take $\psi$ to be the 
inverse of the connection matrix $C$:

\begin{pfms}{\em(of Theorem \ref{thm: duistermaat}).\ \  }
Given $X\in\lp$, let $J:=X/(\pi i)$.
Then we have monodromy data $(b_-,b_+,\Lambda,C):= \wh\nu(1,J)$ as defined in
Section \ref{sn: GW IMs}, (taking $g=1$).
By Lemma \ref{lem: Herm equivariance}, since $J$ is skew-Hermitian,
$C$ is unitary and $b_-=b_+^{-\dagger}$.
The monodromy relation \eqref{eqn: monod reln} implies
\begin{equation}	\label{eqn: new mr}
b^\dagger b = Ce^{2X}C^{-1} = h^\dagger h
\end{equation}
where $b:= b_+$ and $h$ is the Hermitian matrix $Ce^XC^{-1}$.
Now let $h=kan$ be the Iwasawa decomposition of $h$.
Clearly $h^\dagger h=(an)^\dagger an$ and so, from \eqref{eqn: new mr}, we
deduce $b=an$.
Thus 
$$ \wh\delta(Ce^XC^{-1}) = \log(a) = \log(\delta(b)) = \pi i \Lambda,$$
and by definition $\Lambda = \delta(J)$.
Hence if we define $\psi(X)= C^{-1}$ we have established 1).

The real analyticity of $\psi$ is clear: it is the restriction of a
holomorphic map. Property 2) is also straightforward: from Lemma
\ref{lem: inj-etc} we know the map $\wh\nu_J:K\to K^*\times K$ is injective.
Projecting further onto the $K$ factor yields an injective map 
$\pr_K\circ\wh\nu_J:K\to K$ (since the monodromy relation determines the $K^*$
component from the $K$ component). This map is onto for degree reasons and a
diffeomorphism since it is the restriction of a biholomorphic map. 
Finally from 2) of Lemma \ref{lem: inj-etc}, observe that $\phi_X$ is just the
composition of $\pr_K\circ\wh\nu_J$ with the inversion map $K\to K$. \epf
\end{pfms}

Let us briefly continue the story to motivate Theorem \ref{thm: GW ims}.
After Duistermaat, 
the next step was taken by Lu and Ratiu \cite{LuR} who gave $P$ a Poisson
structure by identifying it with the Poisson Lie group $K^*$: The Cartan
decomposition $G=KP$ combined with 
the Iwasawa decomposition $G=KAN$ identifies $P$ with
$AN$, and in turn $K^* \cong AN$. 
Then
the symplectic leaves are the orbits $\cC\subset P$ and 
$\wh \delta$ is a moment map for the dressing action of the maximal
torus $T_K$ of $K$: 
Kostant's non-linear convexity theorem may now be deduced from
the Atiyah, Guillemin and Sternberg convexity theorem.

It was conjectured in \cite{LuR} that there is in fact a 
$T_K$-equivariant Poisson diffeomorphism $\lk^*\cong K^*$.
(So Kostant's non-linear convexity theorem is reduced to the linear case.)
This was proved explicitly for $K=SU(2)$ by P. Xu and then in general by
Ginzburg-Weinstein \cite{GinzW}, building on Duistermaat's indirect 
homotopy argument mentioned above.
Theorem \ref{thm: GW ims} here points out that such diffeomorphisms arise
naturally as monodromy maps for irregular singular connections on the unit 
disc.

\end{section}

\begin{section}{Frobenius Manifolds and Poisson Lie groups}
\label{sn: FMs+PLgps}

Now we will consider the space $U_+$ of Stokes matrices arising in the
theory of Frobenius manifolds.
Our aim is to prove Theorem \ref{thm: fms+plgps} which stated that the
standard Poisson structure on $G^*$ induces the Dubrovin-Ugaglia
Poisson structure on $U_+$.

\begin{pfms}{\em(of Theorem \ref{thm: fms+plgps}).\ \  }
The space $U_+$ appears by restricting the monodromy map to the 
skew-symmetric (complex) matrices, as can be seen from the following:
\begin{lem}	\label{lem: ssym involn}
The monodromy map $\nu$ intertwines the following involutions 
\begin{align*}
i_{\g^*}:\g^*\to&\,\g^*;&  i_{G^*}:G^*\to& \,G^*,\\
B\mapsto -&B^T &
(b_-,b_+,\Lambda) \mapsto& (b_+^T,b_-^T,-\Lambda).
\end{align*}
of $\g^*$ and $G^*$.
In other words: $\nu\circ i_{\g^*} = i_{G^*}\circ \nu$.
\end{lem}
\begin{pf}
This is similar to Lemma \ref{lem: Herm equivariance};
just modify the involutions appearing there to be 
$i(z)=-z, F^\iota := i^*(F^{-T}), A^\iota := - i^*(A^T)$,
and then the proof goes the same:
The original fundamental solutions $\Phi, \Psi$ become 
$\Psi^\iota, \Phi^\iota$. (Using the fact that, under $\iota$, the
function $z^\Lambda$ on $\Sect_0$ becomes 
$(z^\Lambda)^\iota = z^{-\Lambda}e^{-\pi i \Lambda}$ on $\Sect_l$.) 
The lemma is now immediate: for example the relation 
$\Psi b_+ = \Phi$ along the direction $\theta$ implies 
$\Phi^\iota b_+^T = \Psi^\iota$ along $-\theta$.\epf
\end{pf}

Thus the map $U_+\hookrightarrow G^*;\ S\mapsto (S^T,S,0)$ identifies $U_+$
with the fixed point set of the involution $i_{G^*}$ of $G^*$.
(Note that $U_+$ is not embedded as a subgroup.) 
Let $\h\subset \g$ denote the set of skew-symmetric matrices and
identify $\h^*$ with $\h$ using the trace: $ B \leftrightarrow
\tr(B\,\cdot\,)$.
Then Lemma \ref{lem: ssym involn} implies the monodromy map restricts
to a map 
$$\nu\vert_{\h^*}:\h^*\to U_+.$$
This is generically a local analytic isomorphism and the
Dubrovin-Ugaglia Poisson structure on $U_+$ is characterised by the
fact that $\nu\vert_{\h^*}$ is a Poisson map for any value of the irregular
type $A_0$.
(The diagonal entries of $A_0$ are the `canonical coordinates' in the
language of Frobenius manifolds.)
Thus we have a commutative diagram:
$$%
\begin{array}{ccc}
  \g^*  & \mapright{\nu} & G^* \\
\mapup{} &&  \mapup{} \\
\h^* & \mapright{\nu\vert_{\h^*}} & U_+
\end{array}
$$%
where the vertical maps are the inclusions.
By Theorem \ref{thm: poissonness} the top map is Poisson and by
definition the bottom map is Poisson 
(where $\g^*,\h^*$ have their standard complex Poisson structures,
$U_+$ has the Dubrovin-Ugaglia structure and $G^*$ has its standard 
Poisson structure, but scaled by $2\pi i$).  

Now to complete the proof of the theorem we just need to make the
simple observation that
the Poisson structure on $\h^*$ is `induced' from that on $\g^*$ 
via the involution $i_{\g^*}$ (in the sense described after the
statement of the theorem).
\epf
\end{pfms}

For example in the $4\times 4$ case Ugaglia's remarkable 
explicit description \cite{Ugag}
of the Dubrovin-Ugaglia Poisson 
structure on $U_+$ 
is $\{\,\cdot\,,\,\cdot\,\} _{\text{DU}} 
= \frac{\pi i}{2}\{ \,\cdot\,, \,\cdot\,\}$ 
where $\{ \,\cdot\,, \,\cdot\,\}$ 
is determined by the following formulae for its values on
coordinate functions:
$$
\begin{array}{ccc}
S:=\left(\begin{smallmatrix}
1 & u  &  v  & w  \\
0 & 1  &  x  & y  \\
0 & 0  &  1  & z  \\
0 & 0  &  0  & 1 
\end{smallmatrix}\right) &\ \quad\ &
\begin{array}{cl}
\{u,z\}&=\ 0\\
\{v,y\}&=\ 2uz-2xw\\
\{w,x\}&=\ 0
\end{array}
\end{array}
$$
\vspace{-0.6cm}	%
\begin{align*}
\{u,v\}&=2x-uv	&	\{u,w\}&=2y-uw	&\{x,u\}&=2v-xu	&\{y,u\}&=2w-yu	\\
\{v,w\}&=2z-vw	&	\{v,x\}&=2u-vx	&\{z,v\}&=2w-zv	&\{w,y\}&=2u-wy	\\
\{w,z\}&=2v-wz	&	\{x,y\}&=2z-xy	&\{z,x\}&=2y-zx &\{y,z\}&=2x-yz.
\end{align*}
\begin{rmk}
Although our proof is transcendental in nature, the relationship
between the Poisson structures on $G^*$ and $U_+$ is entirely
algebraic.
Indeed, as a plausibility check for Theorem \ref{thm: poissonness} of
this paper, in \cite{thesis}
we used a computer algebra program (Mathematica) to derive
Ugaglia's formula above, from the Poisson structure on $G^*$.
\end{rmk}
\end{section}
\newpage
\renewcommand{\baselinestretch}{1}		%
\normalsize
\bibliographystyle{amsplain}	\label{biby}
\bibliography{../thesis/syr}	
\end{document}